\documentclass[10pt,a4paper]{amsart}
\usepackage[utf8]{inputenc}
\usepackage{hyperref}
\usepackage{amssymb}
\usepackage{bigints}

\newtheorem{theorem}{Theorem}
\newtheorem*{maintheorem*}{Main Theorem}

\def\mgn{{\mathcal M}_{g,n}}

\def\kao{\kappa_1}

\def\wp{Weil-Pe\-ters\-son}

\def\ol#1{{\overline{#1}}}

\newcommand{\cM}{{\mathcal M} }

\newcommand{\wt}{\widetilde}

\def\mgb{{\ol{\cM}_g}}

\newtheorem{remark}{Remark}
\newtheorem{lemma}{Lemma}

\begin{document}
\title[Estimates for the average scalar curvature of the Weil-Petersson metric on $\ol \cM_g$]{Estimates for the average scalar curvature of the Weil-Petersson metric on the Moduli space $\ol \cM_g$}
\author{Georg Schumacher}
\address{Fachbereich Mathematik und Informatik der Philipps-Universit\"at,
Hans-Meerwein-Strasse, Lahnberge, 35032 Marburg, Germany}
\email{schumac@mathematik.uni-marburg.de}
\author{Stefano Trapani}
\address{Dipartimento di Matematica, Universit\'a di Roma, ``Tor Vergata''
Via della Ricerca Scientifica, 00179 Roma, Italy}
\email{trapani@axp.mat.uniroma2.it}

\begin{abstract}
  We give a precise estimate for the average scalar curvature of the Weil-Petersson metric on the moduli space $\ol \cM_g$ as $g\to\infty$ up to the order $1/g^2$.
\end{abstract}

\maketitle

\section{Introduction and Statement of the Result}

The curvature of the \wp\ metric recently attracted further interest. In this note we will give the precise estimate for the average of the scalar curvature $S_{WP}$ of the Weil-Petersson metric on the moduli space $\ol\cM_g$ as $g$ tends to infinity. The result is the value
$$
 \frac{1}{(g-1)}  \frac{\bigintss_\mgb \, (-S_{W\! P})\,dV_{W \! P}}{\bigintss_\mgb dV_{W\! P}} =  \frac{13}{4\pi} +  \frac{\pi}{12} \frac{1}{g} + \left( \frac{1}{4\pi} +\frac{\pi}{12} \right) \frac{1}{g^2}+ O\left(\frac{1}{g^3}\right).
$$
The proof of the asymptotics will be based upon methods of Algebraic Geometry. Wolpert showed in \cite{wo1,wo2} that Mumford's canonical class $\kao$ (the tautological class obtained from the universal curve) from \cite{mu1,mu2} (cf.\ \cite{ac,ac2}) is the cohomology of the \wp\ form extended to the compactification  up to the factor $2\pi^2$ together with the fact that its restriction to the boundary equals to the \wp\ cohomomology of the boundary (interpreted as related to  moduli of punctured surfaces of lower genus).

The finiteness of the \wp\ volume itself is a consequence of Masur's estimates \cite{ma}, whereas the curvature was computed by Wolpert \cite{wo3} and Fischer-Tromba \cite{tro}. These results implied strong negativity properties, in particular the strict negativity of the scalar curvature. It is known that the scalar curvature tends to $-\infty $ towards the boundary. Precise estimates of the curvature of the \wp\ metric towards the boundary are contained in \cite{s} \and \cite{t} with a later developments by Liu-Sun-Yau \cite{lsy1,lsy2}.

Estimates of the \wp\ volume had been given by Mirzakhani \cite{mi}, Mirzakhani-Zograf \cite{m-z}, Penner \cite{pe}, Grushevsky \cite{gru}, Zograf \cite{zo2,zo} and previously in \cite{s-t}. The algebraic aspect is contained in the push-pull formulas by Arbarello and Cornalba \cite{ac,ac2}.

The \wp\ volume of the moduli spaces  $\ol\cM_{g,n}$ of Riemann surfaces of genus $g$ with $n$ punctures is denoted by
$$
V_{g,n}=\int_{\mgn}\kao^{3g-3+n}.
$$
Finally the relationship of intersection numbers and volumes as related to two dimensional gravity ought to pointed out. Pertinent references are \cite{dij,fp,getz,kon,mz,witt}.

We showed the following estimates.
\begin{theorem}[\cite{s-t}]
\begin{itemize}
\item[(i)]
 Let $g>1$. Then
\begin{equation}\label{eq:upper_est}
V_{g,0} \geq \frac{1}{28} V_{g-1,2} + \frac{1}{672} V_{g-1,1} +
\frac{1}{14} \sum_{j=2}^{[g/2]} V_{j,1}V_{g-j,1}
-\frac{1}{28} (V_{\frac{g}{2},1})^2 ,
\end{equation}
with $V_{\frac{g}{2},1}=0$, if $g$ is odd.
\item[(ii)]
There exist constants $0 < c < C$, independent of $n$ such that
\begin{equation}\label{eq:asymp1}
c^g  (2g)! \leq \frac{V_{g,n}}{(3g-3+n)!} \leq C^g  (2g)!
\end{equation}
for all fixed $n\geq 0$ and large $g$.
\end{itemize}
\end{theorem}
Concerning \eqref{eq:asymp1}, a lower estimate for $n=1$ is due to Penner \cite[Theorem 6.2.2]{pe}, and the upper estimates for $n \geq 1$ were first shown by Grushevsky in \cite[Sec.~7]{gru}.

Bounds for the curvature of the \wp\ metric were proven by Wu and Wolf in \cite{ww} and Wu in \cite{wu1,wu2}. A recent result is the following:
\begin{theorem}[{Bridgeman-Wu, \cite{b-w}}]\label{th:briwu}
Denote by $S_{W\!P}$ the scalar curvature of the \wp\ form $\omega_{WP}$, and by $dV_{W\!P}$ its volume element. There exist constants $0<c<C$ such that
\begin{equation}\label{eq:briwu}
c\cdot g \leq \frac{\bigintss_\mgb \, (-S_{W\! P})\,dV_{W \! P}}{\bigintss_\mgb dV_{W\! P}} \leq C\cdot g
\end{equation}
for all $g$.
\end{theorem}
We will show that the Bridgeman-Wu estimate follow in an algebraic way using \eqref{eq:upper_est}.

The actual asymptotics of the \wp\ volume (containing \eqref{eq:asymp1}) was computed by Mirzakhani and Zograf. (We use the above normalization for the \wp\ volume.)
\begin{theorem}[{\cite[Theorem 1.2]{mz}}]
There exists a constant $C\in(0,\infty)$ such that for any given $k\geq 1$, $n\geq 0$
\begin{equation}\label{eq:mz}
  V_{g,n}= C_{MZ} \frac{(3g-3+n)! \, (2g-3+n)! \, 2^{g-3+n}}{\pi^{2g}\sqrt{g}}\left( 1+ \sum_{j=1}^{k} \frac{c_n^{(j)}}{g^j} + O\left(\frac{1}{g^{k+1}}\right)   \right)
\end{equation}
as $g\to\infty$.
\end{theorem}
The theorem contains further characterizations of the polynomials $c^{(j)}_n$ -- we will need the case $k\leq 1$, and $n\leq 2$. The Mirzakhani-Zograf constant $C_{MZ}$ is conjectured to be $C_{MZ}=1/\sqrt\pi$.

Let again $\delta$ be the $\mathbb Q$-divisor related to the boundary components of the moduli space (see also below), and denote by $\kao$ the Mumford class. Recall that the class $[\omega_{WP}]$  of the \wp\ form was computed in \cite{wo1,wo2} as $2 \pi^2 \kao$. Set $\eta= [-Ric(\omega_{WP})/2\pi]$.
\begin{maintheorem*}
  The average total scalar curvature of the \wp\ metric on the moduli space $\ol\cM_g$ satisfies the estimate
  \begin{eqnarray}
   \frac{13}{4\pi} &\leq&
     \frac{3}{\pi} \frac{ \eta \cdot \kao^{3g - 4}}{  \kao^{3g - 3}} = \frac{1}{(g-1)}  \frac{\bigintss_\mgb \, (-S_{W\! P})\,dV_{W \! P}}{\bigintss_\mgb dV_{W\! P}} =  \frac{13}{4\pi} +  \frac{1}{4\pi}  E_g  , \label{eq:Eg}  \quad \text{ where } \\ \nonumber
 \frac{1}{4\pi} E_g &=&  \frac{\pi}{12} \frac{1}{g} + \left( \frac{1}{4\pi} +\frac{\pi}{12} \right) \frac{1}{g^2}+ O\left(\frac{1}{g^3}\right) \text{ for } g\to \infty.\\\text{ The precise value is } \nonumber\\
\label{eq:kappa}
 E_g &=& \frac{\kao^{3g-4} \cdot \delta }{\kao^{3g-3}} \text{ where }\delta = \sum_{j=0}^{[g/2]} \delta_j.
\end{eqnarray}
\end{maintheorem*}
Note that the ampleness of $\kao$ implies that $E_g$ is positive.

\section{Proof}

Let $D$ denote the compactifying divisor of $\cM_g$ with components $\Delta_j$ for $j=0,\ldots,[g/2]$. These give rise to $\mathbb Q$-divisors $\delta_j$ such that in terms of the generally used notation
$$
[\Delta_1]= 2 \delta_1 \text{ and } [\Delta_j]=  \delta_j \text{ for } j \neq 1.
$$
Then $\delta:= \sum_{j=0}^{[g/2]} \delta_j$ so that
$$
D= \delta +\delta_1.
$$
Note that there exist branched two-sheeted coverings $\ol\cM_{g-1,2} \to \Delta_0$, \quad $\ol\cM_{g-1,1}\times \ol\cM_{1,1} \to \Delta_1$, and $\ol\cM_{g/2,1}\times\ol\cM_{g/2,1} \to \Delta_{g/2}$ for even $g$ yielding extra factors $1/2$ in Lemma~\ref{le:le1}.

.

Two classical facts are needed. Mumford's direct image bundle $\lambda$ satisfies \cite{mu1}
$$
\kao = 12 \lambda - \delta,
$$
and by \cite{hm}  the canonical bundle $K_{\ol\cM_{g,0}}$ is equal to
$$
K_{\ol\cM_{g,0}} = 13 \lambda - 2\delta - \delta_1.
$$

In \cite[Corollary 5.5]{t} it was shown by means of a a singular Mumford good hermitian metric that $\eta$ is the Chern class of the dual of the logarithmic tangent bundle $T_{{\overline\cM}_{g,0}}({\log} D)$. Therefore
\begin{equation}\label{eq:ric}
  [\eta] = K_{\ol\cM_{g,0}} + [D].
\end{equation}
Now \eqref{eq:Eg} can be shown:
\begin{gather*}
\frac{ \eta \cdot \kao^{3g - 4}}{  \kao^{3g - 3}} = \frac{ (K_{\ol\cM_{g,0}}+\delta+\delta_1) \cdot \kao^{3g - 4}}{  \kao^{3g - 3}} = \frac{(13 \lambda - \delta )\cdot \kao^{3g - 4}}{  \kao^{3g - 3}} = \frac{(\frac{13}{12} \kao + \frac{1}{12}\delta)\cdot \kao^{3g - 4}}{ \kao^{3g - 3}}
\end{gather*}
The main result consists of the estimate for $E_g$. The special value for $V_{1,1}$ etc.\ are taken from \cite{s-t}.
\begin{lemma}\label{le:le1}
\begin{eqnarray}
  \frac{\kao^{3g-4} \cdot \delta_0 }{\kao^{3g-3}} &=& \frac{1}{2} \frac{V_{g-1,2}}{V_{g,0}}
  \\
  \frac{\kao^{3g-4} \cdot \delta_1 }{\kao^{3g-3}} &=&   \frac{1}{2}  \cdot \frac{V_{1,1}\cdot V_{g-1,1}}{V_{g,0}} = \frac{1}{48} \cdot \frac{ V_{g-1,1}}{V_{g,0}}  \\
  \frac{\kao^{3g-4} \cdot \delta_j }{\kao^{3g-3}} &=& \frac{V_{j,1}\cdot V_{g-j,1}}{V_{g,0}} \qquad \text{ for }\quad 2\leq j\ \leq [(g-1)/2]  \\\label{eq:deltaj}
  \frac{\kao^{3g-4} \cdot \delta_{g/2} }{\kao^{3g-3}} &=&\frac{1}{2}  \frac{(V_{g/2,1})^2}{V_{g,0}}   \text{ \ for \ } g \text{ \ even}.
\end{eqnarray}
\end{lemma}
\begin{remark}
 Lemma~\ref{le:le1}, and \eqref{eq:Eg} together with \eqref{eq:upper_est} imply the Bridgeman-Wu Theorem~\ref{th:briwu}.
\end{remark}

We will apply the Mirzakhani-Zograf estimate \eqref{eq:mz}, and the following values contained in \cite[Remark 1.2]{mz}:
\begin{eqnarray}
  c^{(1)}_0 & = & \frac{7}{12} -\frac{17}{6\pi^2} \label{eq:c0}\\
  c^{(1)}_1 & = & \frac{1}{3} - \frac{5}{6\pi^2} \label{e1:c1}\\
  c^{(1)}_2 & = & \frac{1}{12} + \frac{1}{6\pi^2} .\label{eq:c2}
\end{eqnarray}

\begin{lemma}
  \begin{eqnarray}
    \frac{1}{2}\frac{V_{g-1,2}}{V_{g,0}} &=&
    \frac{\pi^2 }{3g} \left( 1 + \left(\frac{3}{\pi^2} + 1 \right)\frac{1}{g} + O\left( \frac{1}{g^2} \right)     \right)
    \label{eq:g-1_2}\\
    \frac{1}{48}\frac{V_{g-1,1}}{V_{g,0}} &=& \frac{1}{144}\sqrt{\frac{g}{g-1}} \frac{\pi^2}{(g-1)(3g-4)(2g-3)}\left(1 + O\left(\frac{1}{g}\right)\right) = O\left(\frac{1}{g^3}\right)\\
    \frac{V_{j,1}\cdot V_{g-j,1}}{V_{g,0}} &=& \frac{ C_{MZ} }{2}\sqrt{\frac{g}{j(g-j)}} \frac{1}{(3g-3) \left(3g-4\atop 3j-2 \right)(2g-3) \left(2g-4\atop 2j-2  \right) }\left(1 + O\left(\frac{1}{g}\right)\right)\label{eq:Vj1}\\ \nonumber && \text{ for } 2\leq j \leq[(g-1)/2]\\
    \frac{1}{2} \frac{(V_{g/2,1})^2}{V_{g,0}} &=&  \frac{ C_{MZ} }{2}{\frac{1}{\sqrt{g}}} \frac{1}{(3g-3) \left(3g-4\atop (3g-4)/2 \right)(2g-3) \left(2g-4\atop g-2  \right) }\left(1 + O\left(\frac{1}{g}\right) \right)\label{eq:gover2}\\ && \nonumber  \text{ for g even}.
  \end{eqnarray}
\end{lemma}
The computation of \eqref{eq:g-1_2} is the following.
\begin{gather*}
  \frac{1}{2}\frac{V_{g-1,2}}{V_{g,0}} = \frac{1}{2} \frac{(3g-4)!\,(2g-3)!\, 2^{g-2} \pi^{2g}\, \sqrt{g}}{(3g-3)!\,(2g-3)!\,2^{g-3}\,\pi^{2g-2}\,\sqrt{g-1}}\cdot \wt C_g =\\ = \frac{\pi^2}{3g} \left(\frac{g}{g-1} \right)^{3/2} \cdot \wt C_g  = \frac{\pi^2}{3 g } \left( 1 + \frac{3}{2 g} + O\left(\frac{1}{g^2}  \right)   \right) \cdot \wt C_g
\end{gather*}
where
\begin{gather*}
  \wt C_g = \left(1 + c^{(1)}_{2}\frac{1}{g-1} + O\left(\frac{1}{g^2}\right)\right) \Big/ \left(1 + c^{(1)}_{0}\frac{1}{g} + O\left(\frac{1}{g^2}\right)\right) \\
\end{gather*}

The coefficients $c^{(1)}_n$ do not contribute to the remaining terms up to order three in $1/g$.

We need an elementary estimate:

Let $l \leq k \leq n$ be positive integers,  then $(n-k)(k-l) \geq 0,$ i.e. $\frac{n-k+l}{l} \geq \frac{n}{k}.$ So we derive the well known inequality
 \begin{equation} \binom{n}{k} = \prod_{l=1}^{k} \frac{n-k+l}{l}  \geq \left(\frac{n}{k}\right)^k.  \label{binomial} \end{equation}

The function
$$
f(x) := \log \left( \frac{h}{x} \right)^x = x ( \log h - \log x)
$$
attains its maximum value at $x=h/e$ so that its minimum will be taken at either end of a given interval.

The binomial coefficients in \eqref{eq:Vj1} attain asymptotic lower estimates for $2\leq j\leq [g/2]$ at $j=2$ namely
$$
 \left(3g-4\atop 3j-2 \right)\geq \left( \frac{3g-4}{4}\right)^4 \quad \text{and} \quad \left(2g-4\atop 2j-2  \right)\geq (g-2)^2.
$$
Hence, summing up all terms of the type \eqref{eq:Vj1} yields a term of growth order at most $O(1/g^7)$. The last term \eqref{eq:gover2} has a at most a negative exponential growth with respect to $g$ and can also be disregarded.

 This proves the main theorem. \qed

{\it Acknowledgements.} The first named author expresses his thanks to the Department of Mathematics at Tor Vergata for their kind hospitality, and for support by PRIN Real and Complex Manifolds: Topology, Geometry and holomorphic dynamics n.2017JZ2SW5, and by DFG SCHU771/5-1 (Singular Hermitian Metrics/Analytic Theory of Moduli Spaces). The second named author is supported, by PRIN Real and Complex Manifolds: Topology, Geometry and holomorphic dynamics n.2017JZ2SW5, and by MIUR Excellence Department Project awarded to the Department of Mathematics, University of Rome Tor Vergata, CUP E83C18000100006. 

\end{document}